\documentclass[twoside,a4paper,12pt,centertags,reqno]{amsart} 
\usepackage{amsmath,amssymb,verbatim,vmargin}
\usepackage{color}
\usepackage{tikz}
\usepackage{tikz-cd}
\usetikzlibrary{matrix,calc}
\usepackage{hyperref}
\hypersetup{colorlinks,bookmarks=true,linktocpage=true,citecolor=blue,linkcolor=magenta}

\usepackage{eulervm}   

\allowdisplaybreaks 
\usepackage{color}
\usepackage[all]{xy} 

\usepackage{mathtools}
\usepackage{MnSymbol} 

\theoremstyle{plain}
\newtheorem{thm}{Theorem}[section]

\theoremstyle{definition}

\newtheorem{rem}[thm]{Remark}

\usepackage{enumerate}
\usepackage{enumitem}

\newcommand{\bRn}{\mathbb{R}^n}

\newcommand{\pd}{\partial}

\newcommand{\bC}{{\mathbb C}}

\newcommand{\bR}{{\mathbb R}}

\newcommand{\cF}{{\mathcal F}}

\newcommand{\cM}{{\mathcal M}}

\newcommand{\fH}{{\mathbf H}}

\def\barint_#1{\mathchoice
            {\mathop{\vrule width 6pt
height 3 pt depth -2.5pt
                    \kern -9.5pt
\intop \kern -4pt}\nolimits_{#1}}%
            {\mathop{\vrule width 5pt height
3 pt depth -2.6pt
                    \kern -6.5pt
\intop \kern -4pt}\nolimits_{#1}}%
            {\mathop{\vrule width 5pt height
3 pt depth -2.6pt
                    \kern -6pt
\intop \kern -4pt}\nolimits_{#1}}%
            {\mathop{\vrule width 5pt height
3 pt depth -2.6pt
          \kern -6pt \intop \kern -4pt}\nolimits_{#1}}}
          
           \def\bariint_#1{\mathchoice
            {\mathop{\vrule width 15pt
height 3 pt depth -2.5pt
                    \kern -15.8pt
\intop \kern -8pt\intop \kern -4pt}\nolimits_{#1}}%
            {\mathop{\vrule width 9pt height
3 pt depth -2.6pt
                    \kern -10.5pt
\intop \kern -8pt\intop \kern -4pt}\nolimits_{#1}}%
            {\mathop{\vrule width 9pt height
3 pt depth -2.6pt
                    \kern -10pt
\intop \kern -8pt\intop \kern -4pt}\nolimits_{#1}}%
            {\mathop{\vrule width 9pt height
3 pt depth -2.6pt
          \kern -8pt \intop \kern -10pt\intop \kern -4pt}
      \nolimits_{  #1}}}

\def\barintlim_#1{\mathchoice
            {\mathop{\vrule width 6pt
height 3 pt depth -2.5pt
                    \kern -8.8pt
\intop \kern -4pt}\limits_{#1}}%
            {\mathop{\vrule width 5pt height
3 pt depth -2.6pt
                    \kern -6.5pt
\intop \kern -4pt}\limits_{#1}}%
            {\mathop{\vrule width 5pt height
3 pt depth -2.6pt
                    \kern -6pt
\intop \kern -4pt}\limits_{#1}}%
            {\mathop{\vrule width 5pt height
3 pt depth -2.6pt
          \kern -6pt \intop \kern -4pt}\limits_{#1}}}
          
           \def\bariintlim_#1{\mathchoice
            {\mathop{\vrule width 15pt
height 3 pt depth -2.5pt
                    \kern -15.8pt
\intop \kern -8pt\intop \kern -4pt}\limits_{#1}}%
            {\mathop{\vrule width 9pt height
3 pt depth -2.6pt
                    \kern -10.5pt
\intop \kern -8pt\intop \kern -4pt}\limits_{#1}}%
            {\mathop{\vrule width 9pt height
3 pt depth -2.6pt
                    \kern -10pt
\intop \kern -8pt\intop \kern -4pt}\limits_{#1}}%
            {\mathop{\vrule width 9pt height
3 pt depth -2.6pt
          \kern -8pt \intop \kern -10pt\intop \kern -4pt}
      \limits_{  #1}}}
          
\renewcommand{\iint}{\int \kern -8pt\int}       

\numberwithin{equation}{section}
\makeatletter
\@namedef{subjclassname@2020}{\textup{2020} Mathematics Subject Classification}
\makeatother

\title[Commutativity of Maximal Regularity Operators]{The Forward Maximal Regularity Operator Commutes with the Backward}

\author{Yi C. Huang} 
\address{School of Mathematical Sciences, Nanjing Normal University, Nanjing 210023, People's Republic of China}
\email{Yi.Huang.Analysis@gmail.com}
\urladdr{https://orcid.org/0000-0002-1297-7674}

\date{\today} 

\keywords{Abstract evolution equations, Fourier transform, functional calculus.}
\subjclass[2020]{Primary 35K90; Secondary 42A38, 47A60.}  
\thanks{Research of the author is supported by the National NSF grant of China (no. 11801274).}

\begin{document}

\begin{abstract}
Given a generator of a bounded analytic semigroup on a Hilbert space, 
we show that the corresponding forward maximal regularity operator commutes with the backward.
In particular, for self-adjoint generators the images under the two maximal regularity operators have equal unweighted Hilbert space norms.
\end{abstract}

\maketitle


\section{Introduction}

We are concerned with the autonomous evolution equations on $\bR_+=(0,\infty)$:
\begin{equation} \label{e:Evol}
\begin{cases}
\pd_tu+Au=f,\\
 u(0)=0
\end{cases}
\end{equation}
and
\begin{equation} \label{e:Evol'}
\begin{cases}
\pd_tv-Av=f,\\
 v(\infty)=0.
\end{cases}
\end{equation}
Here, the source $f$ lies in the unweighted space $L^2(\bR_+;\fH)$, where $\fH$ is an underlying Hilbert space,
and $-A$ is a densely defined, closed linear operator on $\fH$, 
with domain $D(A)$ and generating a bounded analytic semigroup $\{e^{-zA}: |\arg z|<\delta\}$, $0<\delta<\pi/2$.

The so-called maximal regularity for \eqref{e:Evol} in the Hilbert space $\fH$
(also referred to as de Simon's theorem \cite{dSim64}), and similarly for \eqref{e:Evol'}, can be stated as follows: 

\bigskip 

\begin{quote}
there is some constant $C=C(A,\fH)>0$ such that for all $f\in L^2(\bR_+;\fH)$,
\begin{equation} \label{e:QMR-C}
|||Au|||\leq C|||f|||\quad\text{ and }\quad |||Av|||\leq C|||f|||,
\end{equation}
where $|||\cdot|||=\|\cdot\|_{L^2(\bR_+;\fH)}$.
\end{quote}

\bigskip

\noindent Thus via \eqref{e:Evol} (respectively, \eqref{e:Evol'}) and \eqref{e:QMR-C}, 
the solution $u$ (respectively, $v$) lies in $L^2(\bR_+;D(A))\cap \dot H^1(\bR_+;\fH)$.
The ``maximal regularity" conveyed in \eqref{e:QMR-C} means that $\pd_tu$ and $Au$ 
(respectively, $\pd_tv$ and $Av$) enjoy the same regularity as the source $f$.

Recall that the forward maximal regularity operator $\cM^A_+$ is defined by
\begin{equation} \label{e:M+}
\cM^A_+(f)(t)=\int_0^t Ae^{-(t-s)A}f(s)ds,
\end{equation}
while the backward maximal regularity operator $\cM^A_-$ is given by
\begin{equation} \label{e:M-}
\cM^A_-(f)(t)=\int_t^{\infty} Ae^{-(s-t)A}f(s)ds.
\end{equation}
The operators $\cM^A_\pm$ are associated to \eqref{e:Evol}-\eqref{e:Evol'} as for appropriate $f$, we have
$$Au=\cM^A_+(f)\quad\text{ and }\quad Av=-\cM^A_-(f).$$
Therefore, maximal regularity problems translate into boundedness of $\cM^A_\pm$,
which are typical examples of singular integral operators with operator-valued kernels.

There are so far many machineries proving de Simon's theorem, typically via Fourier transform,
and various extensions to Banach spaces, say, in $L^q(\bR_+;L^p(\Omega))$, where $\Omega\subset\bRn$. 
For more comprehensive materials about the maximal regularity of evolution equations and its applications,
see for example \cite{DenHiePru03, Are04, KunWei04, Mon09} by Denk-Hieber-Pr\"uss, Arendt, Kunstmann-Weis, and Monniaux. 

\section{Commutativity of maximal regularity operators}

In their innovative semigroup approach to generalized Cauchy-Riemann systems and elliptic boundary value problems,
Auscher and Axelsson \cite{AusAxe11} established the endpoint weighted estimates (in $L^2(\bR_+,t^{\mp1}dt;\fH)$) for $\cM^A_\pm$ in presence of quadratic estimates
(see also Hyt\"onen-Ros\'en \cite{HytRos12} for related Carleson duality results, and Ros\'en \cite{Ros12} for the intermediate case with temporal weights $t^{\alpha}dt$ with $|\alpha|<1$). 

In a subsequent proceeding note,
Auscher and Axelsson \cite{AusAxe11p} singled out the abstract formulations of their maximal regularity results obtained in \cite{AusAxe11}.
However, the two maximal regularity operators $\cM^A_\pm$ or rather, the two evolution problems \eqref{e:Evol}-\eqref{e:Evol'}, were treated separately,
which is not the case in \cite{AusAxe11, Ros12}.

The aim of this paper is to consider the ``interaction" of the two maximal regularity operators $\cM^A_\pm$.
Our finding is rather elementary and can be summarized as below.

\begin{thm} \label{thm:commu}
The maximal regularity operators in \eqref{e:M+}-\eqref{e:M-} commute, namely,
\begin{equation} \label{e:commu}
[\cM^A_+,\cM^A_-]=0.
\end{equation}
In particular, if $A=A^*$, we have
\begin{equation} \label{e:equal}
|||\cM^A_+(f)|||=|||\cM^A_-(f)|||
\end{equation}
for all $f\in L^2(\bR_+;\fH)$.
\end{thm}

As a general duality principle, the study of $\cM^A_-$ is reduced to $\cM^{A^*}_+$, see for example \cite{AusAxe11p}.
Here we work directly with $\cM^A_-$ via the traditional Fourier transform approach,
just repeating the arguments (for example \cite[Theorem 2.6]{Mon09}) in proving de Simon's theorem for $\cM^A_+$.
This is indeed the only new ingredient of this paper.

\begin{proof}
Let $\cF$ denote the Fourier transform on $\bR$.
For $\sigma\in\bR$, we have
$$\cF(\cM^A_+(f))(\sigma)=A(i\sigma+A)^{-1}\cF(f)(\sigma),$$ 
and similarly,
$$\cF(\cM^A_-(f))(\sigma)=A(-i\sigma+A)^{-1}\cF(f)(\sigma).$$
The proof of the second identity can be adapted from the arguments for the first (see e.g. \cite[Theorem 2.6]{Mon09}). 
For completeness we give the details.
Let $f\in L^2(\bR_+;\fH)$ and extend $f$ by 0 on $(-\infty,0)$. 
Let $k(t)=Ae^{-tA}$ if $t>0$ and $k(t)=0$ if $t\leq0$.
Thus
$$\cM^A_-(f)(t)=\int_{-\infty}^{\infty}k(s-t)f(s)ds,\quad t\geq0.$$
Moreover, for $\sigma\in\bR$, we have
$$\begin{aligned}
\cF(\cM^A_-(f))(\sigma)&=\int_{-\infty}^{\infty}\int_{-\infty}^{\infty}e^{-it\sigma}k(s-t)f(s)dsdt\\
&=\int_{-\infty}^{\infty}\int_{-\infty}^{\infty}e^{-i(s-t)\sigma}k(t)f(s)dsdt\\
&=\left(\int_{-\infty}^{\infty}e^{it\sigma}k(t)dt\right)\cF(f)(\sigma)\\
&=\left(\int_{0}^{\infty}e^{it\sigma}k(t)dt\right)\cF(f)(\sigma)\\
&=A(-i\sigma+A)^{-1}\cF(f)(\sigma).
\end{aligned}$$
Now, for all $f\in L^2(\bR_+;\fH)$, we compute as follows
\begin{equation} \label{e:comp}
\begin{aligned}
\cF&\left([\cM^A_+,\cM^A_-]f\right)(\sigma)\\
&\quad=\cF(\cM^A_+(\cM^A_-(f)))(\sigma)-\cF(\cM^A_-(\cM^A_+(f)))(\sigma)\\
&\quad=A(i\sigma+A)^{-1}\cF(\cM^A_-(f))(\sigma)-A(-i\sigma+A)^{-1}\cF(\cM^A_+(f))(\sigma)\\
&\quad=A(i\sigma+A)^{-1}A(-i\sigma+A)^{-1}\cF(f)(\sigma)\\
&\qquad\quad-A(-i\sigma+A)^{-1}A(i\sigma+A)^{-1}\cF(f)(\sigma)\\
&\quad=A^2(\sigma^2+A^2)^{-1}\cF(f)(\sigma)-A^2(\sigma^2+A^2)^{-1}\cF(f)(\sigma)=0.\\
\end{aligned}
\end{equation}
Hence, by the isometry via Fourier transform, we have: for all $f\in L^2(\bR_+;\fH)$,
$$[\cM^A_+,\cM^A_-]f=0\quad \text{in}\quad L^2(\bR_+;\fH).$$
We thus proved the (strong) operator identity \eqref{e:commu}.
In view of the adjoint relations
$$(\cM^A_\pm)^*=\cM^{A^*}_\mp,$$
the norm equality \eqref{e:equal} follows from $A=A^*$ and \eqref{e:commu}. 
The theorem is proved.
\end{proof}

\begin{rem}
For the moment we do not know how to apply Theorem \ref{thm:commu} or solely \eqref{e:equal} to the maximal regularity problem on the line or on $\bR_+$
as considered for example in Arendt-Duelli \cite{AreDue06}, Arendt-Zamboni \cite{AreZam10}, and Auscher-Axelsson \cite{AusAxe11}.
In particular, we expect our elementary finding in Theorem \ref{thm:commu} can be connected to the Cauchy non-integral formulae systematically studied in Ros\'en \cite{Ros12}.
\end{rem}

\section{Commutativity of holomorphic functional calculus}

Conceptually one may think of the commutativity formula \eqref{e:commu} as an algebra-homomorphism property of Albrecht's operational calculus (recalled in \cite{AusAxe11}).
Using Fourier transform, we justified this intuition at the level of functional calculus.

Let us take a further look at the calculus for the Fourier multipliers: 
$$m_\pm(\sigma,\xi)=\xi(\pm i\sigma+\xi)^{-1}\quad\longmapsto \quad m_\pm(\sigma,A)=A(\pm i\sigma+A)^{-1}.$$ 
So what we did in above proof of Theorem \ref{thm:commu}, and in particular in the computations \eqref{e:comp}, 
is actually to expose the rather elementary commutativity
\begin{equation} \label{e:commu'}
[m_+(\sigma,A),m_-(\sigma,A)]=0.
\end{equation}
Note that the two functions $m_\pm(\sigma,\xi)$, parametrized by $\sigma\in\bR$ and holomorphic on $\bC\backslash i\bR$ , have no decay when $|\xi|\rightarrow\infty$.
Therefore, the commutativity \eqref{e:commu'} does \textit{not} belong to the algebra-homomorphism property of holomorphic calculus 
(see e.g. \cite[Proposition 6.2]{AusAxe11} in the operational calculus context).
However, under McIntosh's bounded holomorphic functional calculus for $A$, 
we do have such a property which roughly states that for two bounded holomorphic functions $b_1$ and $b_2$, 
$$b_1(A)b_2(A)=(b_1b_2)(A),$$
see Albrecht-Duong-McIntosh \cite{AlbDuoMcI96} and Haase \cite{Haas06}.
Since $b_1b_2=b_2b_1$, we have $[b_1(A),b_2(A)]=0$ as bounded operators acting on $\fH$. Thus we also have
\begin{equation} \label{e:commu''}
[\cM^A_+b_1(A),\cM^A_-b_2(A)]=0
\end{equation}
as bounded operators acting on $L^2(\bR_+;\fH)$.
\eqref{e:commu''} does include \eqref{e:commu} as a special case.

In the semigroup approach to Cauchy-Riemann systems and elliptic boundary value problems, 
the candidates $b_1$ and $b_2$ can be the spectral projections of (bisectorial) multiplicatively perturbed Dirac operators.
Note that the non-trivial case happens when $b_1$ and $b_2$ are the same projection.
See \cite{AusAxe11, Ros12} for details.

\bigskip

\section*{\textbf{Compliance with ethical standards}}

\bigskip

\textbf{Conflict of interest} The author has no known competing financial interests or
personal relationships that could have appeared to influence this reported work.

\bigskip

\textbf{Availability of data and material} Not applicable.

\bigskip

\bibliographystyle{alpha}
 
\bibliography{Hua-MaxRegCommute}

\end{document}